\documentclass[11pt]{article}
\usepackage{amsfonts}
\usepackage{latexsym}

\newcommand{\qed}{\hfill$\Box$}

\newcommand{\then}{\Longrightarrow}
\newcommand{\inv}{{\it inv}}
\newcommand{\bpi}{{\bar\pi}}

\newtheorem{thm}{Theorem}[section]
\newtheorem{pro}[thm]{Proposition}
\newtheorem{lem}[thm]{Lemma}

\newtheorem{cor}[thm]{Corollary}

\newtheorem{obs}[thm]{Observation}

\newtheorem{exa}[thm]{Example}
\newtheorem{df}[thm]{Definition}
\newtheorem{rem}[thm]{Remark}

\begin{document}
\pagestyle{myheadings}

\title{On Degrees in the Hasse Diagram\\ of the Strong Bruhat Order}
\bibliographystyle{acm}
\author{Ron M. Adin%
\thanks{Department of Mathematics and Statistics, Bar-Ilan University,
Ramat-Gan 52900, Israel. Email: {\tt radin@math.biu.ac.il}}\ $^\S$
\and Yuval Roichman%
\thanks{Department of Mathematics and Statistics, Bar-Ilan University,
Ramat-Gan 52900, Israel. Email: {\tt yuvalr@math.biu.ac.il}}
\thanks{Research of both authors was supported in part by
the Israel Science Foundation, founded by the Israel Academy of
Sciences and Humanities, and by the EC's IHRP Programme, within the
Research Training Network ``Algebraic Combinatorics in Europe'',
grant HPRN-CT-2001-00272.}}
\date{March 12, 2006}

\maketitle

\begin{abstract}
For a permutation $\pi$ in the symmetric group $S_n$ let the
%{\it strong descent number}, $d_{-}(\pi)$,
{\it total degree} %, $d(\pi)$,
be its valency in the Hasse diagram of the strong Bruhat order on $S_n$,
and let the {\it down degree} %, $d_{-}(\pi)$,
be the number of permutations which are covered by $\pi$ in 
the strong Bruhat order.
The maxima of
%the strong descent number and
the total degree and the down degree and their values at a random permutation
are computed. Proofs involve variants of a classical theorem of Tur\'an from
extremal graph theory.
\end{abstract}

\section{The %Strong Descent Number and the
Down, Up and Total Degrees}

%\subsection{Definitions}

\begin{df}
For a permutation $\pi\in S_n$ let the
%{\em strong descent number} (or the
{\em down degree} $d_{-}(\pi)$ be the number of permutations in $S_n$ which are
covered by $\pi$ in the strong Bruhat order. Let the {\em up degree}
$d_+(\pi)$ be the number of permutations which cover $\pi$ in this order.
The {\em total degree} of $\pi$ is the sum
$$
d(\pi):=d_-(\pi)+d_+(\pi),
$$
i.e., the valency of $\pi$ in the Hasse diagram of the strong
Bruhat order.
\end{df}

Explicitly, for $1\le a < b \le n$ let $t_{a,b}=t_{b,a}\in S_n$ be 
the transposition interchanging $a$ and $b$, 
and for $\pi\in S_n$ let
$$
\ell(\pi) := \min\{k \,|\, \pi = 
                   s_{i_1}s_{i_2} \cdots s_{i_k}\}
$$
be the {\em length} of $\pi$ with respect to the standard Coxeter generators
$s_i = t_{i,i+1}$ $(1\le i < n)$ of $S_n$. Then
\begin{eqnarray*}
d_{-}(\pi) &=& \#\{t_{a,b} \,|\, \ell(t_{a,b}\pi) = \ell(\pi) - 1\}\cr
d_{+}(\pi) &=& \#\{t_{a,b} \,|\, \ell(t_{a,b}\pi) = \ell(\pi) + 1\}\cr
d(\pi) = d_{-}(\pi) + d_{+}(\pi) 
           &=& \#\{t_{a,b} \,|\, \ell(t_{a,b}\pi) = \ell(\pi) \pm 1\}
\end{eqnarray*}
For the general definitions and other properties of the weak and strong
Bruhat orders see, e.g., \cite[Ex. 3.75]{ECI} and~\cite[\S\S2.1, 3.1]{BB}.

We shall describe $\pi\in S_n$ by its sequence of values 
$[\pi(1),\ldots,\pi(n)]$.

\begin{obs}\label{t.ob1}
%The permutation $\sigma$ is covered by $\pi$ in the strong Bruhat
%order on $S_n$ if and only if there exist $1\le a< b\le n$ and a
%reflection $t_{a,b}$ such that $\sigma=t_{a,b}\pi$, where
%$\pi=[\dots,b,\dots,a,\dots]$ and there is no letter between the
%positions of $b$ and $a$ in $\pi$ whose value is between $a$ and $b$.
$\pi$ covers $\sigma$ in the strong Bruhat order on $S_n$ if and only if
there exist $1\le i < k \le n$ such that
\begin{enumerate}
\item
$b := \pi(i) > \pi(k) =: a$.
\item
$\sigma = t_{a,b}\pi$, i.e.,
$\pi = [\ldots,b,\ldots,a,\ldots]$ and
$\sigma = [\ldots,a,\ldots,b,\ldots]$.
\item
There is no $i < j < k$ such that $a < \pi(j) < b$.
\end{enumerate}
\end{obs}

\begin{cor}
For every $\pi\in S_n$
$$
d_{-}(\pi)=d_{-}(\pi^{-1}).
$$
\end{cor}

\begin{exa}
In $S_3$,
$d_-[123]=0$, 
$d_-[132]=d_-[213]=1$, and
$d_{-}[321]=d_-[231]=d_-[312]=2$.
On the other hand,
$d[321]=d[123]=2$ and
$d[213]=d[132]=d[312]=d[231]=3$.
\end{exa}

\begin{rem}
The classical descent number of a permutation $\pi$ in the symmetric group
$S_n$ is the number of permutations in $S_n$ which are covered by $\pi$ in the
(right) {\em weak} Bruhat order. Thus, the down degree may be considered as
a ``strong descent number''.
\end{rem}

\begin{df}
For $\pi\in S_n$ denote
$$
D_{-}(\pi) := \{t_{a,b} \,|\, \ell(t_{a,b}\pi)=\ell(\pi)-1\},
$$
the {\it strong descent set} of $\pi$.
\end{df}

\begin{exa}
The strong descent set of $\pi=[7,9,5,2,3,8,4,1,6]$
is
$$
D_{-}(\pi) =
\{t_{1,2}, t_{1,3}, t_{1,4}, t_{2,5}, t_{3,5}, t_{4,5}, 
  t_{4,8}, t_{5,7}, t_{5,9}, t_{6,7}, t_{6,8}, t_{8,9}\}.
$$
\end{exa}

\begin{rem}
Generalized pattern avoidance, involving strong descent sets, was applied by
Woo and Yong~\cite{Yong} to determine which Schubert varieties are Gorenstein.
\end{rem}

\begin{pro}\label{t.determine}
The strong descent set $D_{-}(\pi)$ uniquely determines the permutation $\pi$.
\end{pro}

\noindent{\bf Proof.} 
By induction on $n$. The claim clearly holds for $n=1$.

Let $\pi$ be a permutation in $S_n$, and let $\bpi\in S_{n-1}$ be 
the permutation obtained by deleting the value $n$ from $\pi$. 
Note that, by Observation~\ref{t.ob1},
$$
D_{-}(\bpi)=D_{-}(\pi)\setminus \{t_{a,n} \,|\, 1\le a < n\}.
$$
By the induction hypothesis $\bpi$ is uniquely determined by this set. 
Hence it suffices to determine the position of $n$ in $\pi$.

Now, if $j:= \pi^{-1}(n) < n$ then clearly $t_{\pi(j+1),n}\in D_{-}(\pi)$. 
Moreover, by Observation~\ref{t.ob1}, 
$t_{a,n}\in D_{-}(\pi) \then a\ge\pi(j+1)$. Thus $D_{-}(\pi)$ determines
$$
\bpi(j) = \pi(j+1) = \min\{a \,|\, t_{a,n}\in D_{-}(\pi)\},
$$ 
and therefore determines $j$. 
Note that this set of $a$'s is empty if and only if $j=n$.
This completes the proof.

\qed

\section{Maximal Down Degree}

In this section we compute the maximal value of the
down degree on $S_n$ and find all the permutations achieving the maximum. 
We prove

\begin{pro}\label{t.max}
For every positive integer $n$
$$
\max\{d_{-}(\pi)|\ \pi\in S_n\}= \lfloor n^2/4 \rfloor .
$$
\end{pro}

\begin{rem}
The same number appears as the order dimension of the strong Bruhat 
poset~\cite{Reading}.
An upper bound on the maximal down degree for finite Coxeter groups appears in
\cite[Prop. 3.4]{Br}.
\end{rem}

For the proof of Proposition~\ref{t.max} we need a classical theorem of Tur\'an.

\begin{df}
Let $r\le n$ be positive integers. The {\em Tur\'an graph}
$T_r(n)$ is the complete $r$-partite graph with $n$ vertices and
all parts as equal in size as possible, i.e., each size is 
either $\lfloor n/r \rfloor$ or $\lceil n/r \rceil$. 
Denote by $t_r(n)$ the number of edges of $T_r(n)$.
\end{df}

\begin{thm}\label{t.Turan}{\rm~\cite{Turan}~\cite[IV, Theorem 8]{Bol}
{\bf (Tur\'an's Theorem)}}
\begin{itemize}
\item[(1)] 
Every graph of order $n$ with more than $t_r(n)$ edges
contains a complete subgraph of order $r+1$. 
\item[(2)] 
$T_r(n)$ is
the unique graph of order $n$ with $t_r(n)$ edges that does not
contain a complete subgraph of order $r+1$.
\end{itemize}
\end{thm}

We shall apply 
the special case $r=2$ (due to Mantel) of
Tur\'an's theorem to the following graph.

\begin{df}
The {\em strong descent graph} of $\pi\in S_n$, denoted $\Gamma_{-}(\pi)$,
is the undirected graph whose set of vertices is $\{1,\dots,n\}$ and whose
set of edges is
$$
 \{\{a,b\} \,|\, t_{a,b}\in D_{-}(\pi)\}.
$$
\end{df}

By definition, the number of edges in $\Gamma_{-}(\pi)$ equals $d_{-}(\pi)$.

\begin{rem}
Permutations for which the strong descent graph is connected are
called {\em indecomposable}. Their enumeration was studied
in~\cite{Comtet}; see~\cite[pp. 7--8]{GS}. 
The number of components in $\Gamma_{-}(\pi)$ is equal to the number of 
{\em global descents} in $\pi w_0$ (where $w_0:=[n,n-1,\dots,1]$), 
which were introduced and studied in~\cite[Corollaries 6.3 and 6.4]{Sottile}.
\end{rem}

%\medskip

\begin{lem}\label{t.triangle}
For every $\pi\in S_n$, the strong descent graph $\Gamma_{-}(\pi)$ is triangle-free.
\end{lem}

\noindent{\bf Proof.} 
Assume that $\Gamma_{-}(\pi)$ contains a triangle. 
Then there exist $1\le a < b < c\le n$ 
such that $t_{a,b}, t_{a,c}, t_{b,c} \in D_{-}(\pi)$. %are strong descents of $\pi$. 
By Observation~\ref{t.ob1}, 
$$
t_{a,b}, t_{b,c} \in D_{-}(\pi) \then
\pi^{-1}(c) < \pi^{-1}(b) < \pi^{-1}(a) \then
t_{a,c} \not\in D_{-}(\pi).
$$
This is a contradiction.

\qed

\medskip

\noindent{\bf Proof of Proposition~\ref{t.max}.} 
By Theorem~\ref{t.Turan}(1) together with Lemma~\ref{t.triangle}, for every
$\pi\in S_n$
$$
d_{-}(\pi)\le t_2(n)=\lfloor n^2/4 \rfloor .
$$
Equality holds since
$$
d_{-}([\lfloor n/2 \rfloor+1,\lfloor n/2 \rfloor+2,\dots,n,
1,2,\dots,\lfloor n/2 \rfloor]) =
\lfloor n^2/4 \rfloor .
$$

\qed

\medskip

Next we classify (and enumerate) the permutations which achieve
the maximal down degree. %strong descent number.

\begin{lem}
Let $\pi\in S_n$ be a permutation with maximal down degree. %strong descent number. 
Then $\pi$ has no decreasing subsequence of length 4.
\end{lem}

\noindent{\bf Proof.}
Assume that $\pi = [\ldots d \ldots c \ldots b \ldots a \ldots]$ 
with $d>c>b>a$ and $\pi^{-1}(a) - \pi^{-1}(d)$ minimal. Then 
%$t_{a,b},t_{b,c},$ and $t_{c,d}$ are strong descents of $\pi$.
$t_{a,b},t_{b,c},t_{c,d} \in D_{-}(\pi)$ but,
by Observation~\ref{t.ob1}, $t_{a,d} \not\in D_{-}(\pi)$.
It follows that $\Gamma_{-}(\pi)$ is not a complete bipartite graph,
since $\{a,b\}$, $\{b,c\}$, and $\{c,d\}$ are edges but $\{a,d\}$ is not. 
By Lemma~\ref{t.triangle}, combined with Theorem~\ref{t.Turan}(2),
the number of edges in $\Gamma_{-}(\pi)$ is less than $\lfloor n^2/4\rfloor$.

\qed

\smallskip

\begin{pro}
For every positive integer $n$
$$
\#\{\pi\in S_n \,|\, d_{-}(\pi) = \lfloor n^2/4 \rfloor\} =
\cases{n,& if $n$ is odd;\cr
n/2,& if $n$ is even.}
$$
Each such permutation has the form
$$
\pi = [t+m+1, t+m+2, \ldots, n, t+1, t+2, \ldots, t+m, 1, 2, \ldots, t],
$$
where $m\in \{\lfloor n/2\rfloor, \lceil n/2 \rceil\}$ and $1\le t \le n-m$.
Note that $t=n-m$ (for $m$) gives the same permutation as $t=0$ (for $n-m$
instead of $m$).
\end{pro}

\noindent{\bf Proof.}
It is easy to verify the claim for $n\le 3$. Assume $n\ge 4$.

Let $\pi\in S_n$ with $d_{-}(\pi)=\lfloor n^2/4\rfloor$. 
By Theorem~\ref{t.Turan}(2), $\Gamma_{-}(\pi)$ is isomorphic to 
the complete bipartite graph  $K_{\lfloor n/2\rfloor, \lceil n/2 \rceil}$. 
Since $n\ge 4$, each side of the graph contains at least two vertices. 
Let $1 = a < b$ be two vertices on one side, and $c < d$ two vertices on the
other side of the graph.
%that contains the vertex $1$, and $c,d$ two of the vertices in the
%other side. The subgraph induced on $\{1,b,c,d\}$ is isomorphic to
%$K_{2,2}$ where $1,b$ lie in one side and $c,d$ in the other.
%W.l.o.g. $d>c$. We have three cases  to check:
Since $t_{b,c}, t_{b,d}\in D_{-}(\pi)$, there are three possible cases:

\begin{enumerate}
\item
$b<c$, and then $\pi = [\ldots c \ldots d \ldots b \ldots]$
(since $\pi = [\ldots d \ldots c \ldots b \ldots]$ contradicts 
$t_{b,d}\in D_{-}(\pi)$).
\item
$c<b<d$, and then $\pi = [\ldots d \ldots b \ldots c \ldots]$.
\item
$d<b$, and then $\pi = [\ldots b \ldots c \ldots d \ldots]$
(since $\pi = [\ldots b \ldots d \ldots c \ldots]$ contradicts 
$t_{b,c}\in D_{-}(\pi)$).
\end{enumerate}

The same also holds for $a=1$ instead of $b$, but then 
cases $2$ and $3$ are impossible since $a=1<c$.
Thus necessarily $c$ appears before $d$ in $\pi$, and case $2$ is therefore 
impossible for {\em any} $b$ on the same side as $a=1$.
In other words: no vertex on the same side as $a=1$ is intermediate, 
either in position (in $\pi$) or in value, to $c$ and $d$.

Assume now that $n$ is even.
The vertices not on the side of $1$ form (in $\pi$) a block of length $n/2$
of numbers which are consecutive in value as well in position. 
They also form an increasing subsequence of $\pi$, since $\Gamma_{-}(\pi)$ 
is bipartite. 
The numbers preceding them are all larger in value, and are increasing; 
the numbers succeeding them are all smaller in value, are increasing, 
and contain $1$. It is easy to check that each permutation $\pi$ of this form
has maximal $d_{-}(\pi)$. Finally, $\pi$ is completely determined by
the length $1\le t \le n/2$ of the last increasing subsequence.

For $n$ odd one obtains a similar classification, except that the length
of the side not containing $1$ is either $\lfloor n/2 \rfloor$ or
$\lceil n/2 \rceil$. This completes the proof.

\qed

\section{Maximal Total Degree}

%Proposition~\ref{t.max} also holds, of course, for $d_{+}$ instead of $d_{-}$.
%It follows that 
Obviously, the maximal value of the total degree $d = d_{-} + d_{+}$
cannot exceed ${n \choose 2}$, the total number of transpositions in $S_n$.
This is slightly better than the bound $2\lfloor n^2/4 \rfloor$
obtainable from Proposition~\ref{t.max}. 
The actual maximal value is smaller.
 
\begin{thm}\label{t.total} 
For $n\ge 2$,
the maximal total degree in the Hasse diagram of the
strong Bruhat order on $S_n$ is
$$
\lfloor n^2/4 \rfloor +n-2.
$$
\end{thm}

In order to prove this result, associate with each permutation $\pi\in S_n$
a graph $\Gamma(\pi)$, whose set of vertices is $\{1,\dots,n\}$ and
whose  set of edges is
$$
\{\{a,b\} \,|\, \ell(t_{a,b}\pi) -\ell(\pi)=\pm 1\}.
$$
This graph has many properties; e.g., it is $K_5$-free and is the 
edge-disjoint union of two triangle-free graphs on the same set of vertices.
However, these properties are not strong enough to imply the above result.
A property which does imply it is the following bound on the minimal degree.

\begin{lem}\label{t.mindeg}
There exists a vertex in $\Gamma(\pi)$ with degree at most $\lfloor n/2\rfloor+1$. 
\end{lem}

\noindent{\bf Proof.} 
Assume, on the contrary, that each vertex in $\Gamma(\pi)$ has at least
$\lfloor n/2\rfloor +2$ neighbors. This applies, in particular, to the vertex 
$\pi(1)$.
Being the first value of $\pi$, the neighborhood of $\pi(1)$ in $\Gamma(\pi)$,
viewed as a subsequence of $[\pi(2),\ldots,\pi(n)]$,
consists of a shuffle of a decreasing sequence of numbers larger than $\pi(1)$
and an increasing sequence of numbers smaller than $\pi(1)$. 
Let $a$ be the rightmost neighbor of $\pi(1)$. 
The intersection of the neighborhood of $a$ with
the neighborhood of $\pi(1)$ is of cardinality at most two. 
Thus the degree of $a$ is at most
$$
n-(\lfloor n/2\rfloor +2)+2=\lceil n/2\rceil \le \lfloor n/2 \rfloor +1,
$$
which is a contradiction.

\qed

\medskip

\noindent{\bf Proof of Theorem~\ref{t.total}.} 
First note that, by definition, the total degree of $\pi\in S_n$ in the Hasse diagram
of the strong Bruhat order is equal to the number of edges in $\Gamma(\pi)$.
We will prove that this number $e(\Gamma(\pi)) \le \lfloor n^2/4 \rfloor +n-2$, 
by induction on $n$.

The claim is clearly true for $n = 2$. 
Assume that the claim holds for $n-1$, and let $\pi\in S_n$. 
Let $a$ be a vertex of $\Gamma(\pi)$ with minimal degree, 
and let $\bpi\in S_{n-1}$ be the permutation obtained from $\pi$ by 
deleting the value $a$
(and decreasing by $1$ all the values larger than $a$). Then
$$
e(\Gamma(\bpi)) \ge e(\Gamma(\pi)\setminus a),
$$
where the latter is the number of edges in $\Gamma(\pi)$ which are not 
incident with the vertex $a$.
By the induction hypothesis and Lemma~\ref{t.mindeg},
\begin{eqnarray*}
e(\Gamma(\pi)) 
&=& e(\Gamma(\pi)\setminus a) + d(a) \le e(\Gamma(\bpi)) + d(a)\cr 
&\le& \lfloor (n-1)^2/4 \rfloor +(n-1)-2 + \lfloor n/2\rfloor +1\cr
&=& \lfloor n^2/4 \rfloor +n-2.
\end{eqnarray*}
Equality holds since, letting $m := \lfloor n/2 \rfloor$,
$$
e(\Gamma([m+1, m+2, \ldots, n, 1, 2, \dots, m])) = \lfloor n^2/4 \rfloor +n-2.
$$

\qed

\begin{thm}
$$
\#\{\pi\in S_n \,|\, d(\pi) = \lfloor n^2/4 \rfloor +n-2\} = \cases{%
2,& if $n=2$;\cr
4,& if $n=3$ or $n=4$;\cr
8,& if $n\ge 6$ is even;\cr
16,& if $n\ge 5$ is odd.}
$$
The extremal permutations have one of the following forms:
$$
\pi_0 := [m+1, m+2, \ldots, n, 1, 2, \ldots, m]\qquad
(m \in \{\lfloor n/2 \rfloor, \lceil n/2 \rceil\}),
$$
and the permutations obtained from $\pi_0$ by one or more of the following operations:
\begin{eqnarray*}
\pi &\mapsto& \pi^r := [\pi(n), \pi(n-1), \ldots, \pi(2), \pi(1)] 
\qquad\mbox{\rm(reversing $\pi$),}\cr
\pi &\mapsto& \pi^s := \pi \cdot t_{1,n} 
\qquad\mbox{\rm(interchanging $\pi(1)$ and $\pi(n)$),}\cr
\pi &\mapsto& \pi^t := t_{1,n} \cdot \pi 
\qquad\mbox{\rm(interchanging $1$ and $n$ in $\pi$).}
\end{eqnarray*}
\end{thm}

\noindent{\bf Proof.}
It is not difficult to see that all the specified permutations are indeed
extremal, and their number is as claimed (for all $n\ge 2$). 

The claim that there are no other extremal permutations will be proved
by induction on $n$.
For small values of $n$ (say $n\le 4$) this may be verified directly.
Assume now that the claim holds for some $n\ge 4$, 
and let $\pi\in S_{n+1}$ be extremal.
Following the proof of Lemma~\ref{t.mindeg},
let $a$ be a vertex of $\Gamma(\pi)$ with degree at most 
$\lfloor (n+1)/2 \rfloor +1$, which is either $\pi(1)$ or its rightmost neighbor.
As in the proof of Theorem~\ref{t.total},
let $\bpi\in S_n$ be the permutation obtained from $\pi$ by 
deleting the value $a$ (and decreasing by $1$ all the values larger than $a$).
All the inequalities in the proof of Theorem~\ref{t.total} must hold as 
equalities, namely:
$e(\Gamma(\pi)\setminus a)=e(\Gamma(\bpi))$,
$d(a) = \lfloor (n+1)/2 \rfloor +1$,
and $\bpi$ is extremal in $S_n$. 
By the induction hypothesis, $\bpi$ must have one of the prescribed forms.
In all of them, $\{\bpi(1),\bpi(n)\} = \{m,m+1\}$ is an edge of $\Gamma(\bpi)$.
Therefore the corresponding edge $\{\pi(1),\pi(n+1)\}$
(or $\{\pi(2),\pi(n+1)\}$ if $a=\pi(1)$, or $\{\pi(1),\pi(n)\}$ if $a=\pi(n+1)$)
is an edge of $\Gamma(\pi)\setminus a$, namely of $\Gamma(\pi)$.
If $a\ne\pi(1),\pi(n+1)$ then $\pi(n+1)$ is the rightmost neighbor of $\pi(1)$,
contradicting the choice of $a$.
If $a=\pi(n+1)$ we may use the operation $\pi \mapsto \pi^r$.
Thus we may assume from now on that $a=\pi(1)$.

Let $N(a)$ denote the set of neighbors of $a$ in $\Gamma(\pi)$.
Assume first that
$$
\bpi = \pi_0 = [m+1, m+2, \ldots, n, 1, 2, \ldots, m]\qquad
(m \in \{\lfloor n/2 \rfloor, \lceil n/2 \rceil\}).
$$
Noting that $\lceil n/2 \rceil = \lfloor (n+1)/2 \rfloor$ and 
keeping in mind the decrease in certain values during the transition 
$\pi\mapsto\bpi$, we have the following cases:
\begin{description}
\item{(1)} $a > m+1:$
in this case $1,\ldots, m \not\in N(a)$, so that 
$$
d(a) \le n-m \le \lceil n/2 \rceil = \lfloor (n+1)/2 \rfloor
< \lfloor (n+1)/2 \rfloor +1.
$$
Thus $\pi$ is not extremal.
\item{(2)} $a < m:$
in this case $m+3,\ldots, n+1, m+1 \not\in N(a)$, so that
$$
d(a) \le 1 + (m-1) \le \lceil n/2 \rceil < \lfloor (n+1)/2 \rfloor +1.
$$
Again, $\pi$ is not extremal.
\item{(3)} $a\in\{m,m+1\}:$
in this case 
$$
d(a) = 1 + m \le \lfloor (n+1)/2 \rfloor +1,
$$
with equality iff $m = \lfloor (n+1)/2 \rfloor$.
This gives $\pi\in S_{n+1}$ of the required form (either $\pi_0$ or $\pi_0^s$).
\end{description}

A similar analysis for $\bpi = \pi_0^s$ gives extremal permutations
only for $a\in\{m+1,m+2\}$ and $d(a)=3$, so that $n=4$ and $\bpi=[2413]\in S_4$.
The permutations obtained are $\pi = [32514]$ and $\pi = [42513]$, which are 
$\pi_0^{rt}, \pi_0^{rst}\in S_5$, respectively.

The other possible values of $\bpi$ are obtained by the $\pi \mapsto \pi^r$ 
and $\pi \mapsto \pi^t$ operations from the ones above, 
and yield analogous results.

\qed

\section{Expectation}

In this subsection we prove an exact formula for the expectation of 
the down degree of a permutation in $S_n$.

\begin{thm}\label{t.e11} 
For every positive integer $n$, the expected down degree %strong descent number
of a random permutation in $S_n$ is
$$
E_{\pi\in S_n} [d_{-}(\pi)]
= \sum_{i=2}^n \sum_{j=2}^i \sum_{k=2}^j \frac{1}{i\cdot(k-1)}
= (n+1)\sum_{i=1}^n \frac{1}{i} - 2n.
$$
\end{thm}

It follows that

\begin{cor}\label{t.e1}
As $n\to\infty$,
$$
E_{\pi\in S_n} [d_{-}(\pi)] = n \ln n + O(n)
$$
and
$$
E_{\pi\in S_n} [d(\pi)] = 2n \ln n + O(n).
$$
\end{cor}

To prove Theorem~\ref{t.e11} we need some notation.
For $\pi\in S_n$ and $2\le i\le n$ let $\pi_{|i}$ be the permutation 
obtained from $\pi$ by omitting all letters which are larger than 
or equal to $i$.
For example, if $\pi=[6,1,4,8,3,2,5,9,7]$ then
$\pi_{|9}=[6,1,4,8,3,2,5,7]$, $\pi_{|7}= [6,1,4,3,2,5]$, and
$\pi_{|4}=[1,3,2]$.

Also, denote by $\pi^{|j}$ the suffix of length $j$ of $\pi$. 
For example, if $\pi=[6,1,4,8,3,2,5,9,7]$ then
$\pi^{|3}=[5,9,7]$ and $\pi_{|4}^{|2}=[3,2]$.

Let $l.t.r.m.(\pi)$ be the number of left-to-right maxima in $\pi$:
$$
l.t.r.m.(\pi) := \#\{i \,|\, \pi(i) = \max_{1\le j\le i} \pi(j)\}
$$

\begin{lem}\label{t.e12}
For every $\pi\in S_n$, if $\pi_{|{i+1}}^{-1}(i)=j$ then
$$
d_{-}(\pi_{|i+1})-d_{-}(\pi_{|i})= l.t.r.m.(\pi_{|i}^{|i-j}).
$$
\end{lem}

%\noindent
%{\bf Proof.}
%Both sides are equal to 
%$$
%\#\{t \,|\, t<i,  
%$$
%\qed

\noindent{\bf Proof of Theorem~\ref{t.e11}.}
%By Observation~\ref{t.ob1}
Clearly, for every $\pi\in S_n$
$$
d_{-}(\pi)=\sum_{i=2}^{n} \left[d_{-}(\pi_{|i+1})-d_{-}(\pi_{|i})\right].
$$
Thus, by Lemma~\ref{t.e12},
$$%\begin{equation}\label{e.1}
d_{-}(\pi)=\sum_{i=2}^{n} l.t.r.m.(\pi_{|i}^{|i-j_i}),
$$%\end{equation}
where $j_i$ is the position of $i$ in $\pi_{|{i+1}}$, i.e.,
$j_i:=\pi_{|{i+1}}^{-1}(i)$.

Define a random variable $X$ to be the down degree $d_{-}(\pi)$ 
of a random (uniformly distributed) permutation $\pi\in S_n$. 
Then, for each $2\le i\le n$, 
$\pi_{|i+1}$ is a random (uniformly distributed) permutation in $S_i$, 
and therefore 
$j=\pi_{|i+1}^{-1}(i)$ is uniformly distributed in $\{1,\ldots,i\}$
and
$\pi_{|i+1}^{|i-j}$ is essentially a random (uniformly distributed) 
permutation in $S_{i-j}$ (after monotonically renaming its values).
Therefore, by linearity of the expectation,
\begin{equation}\label{e.2}
E[X] = \sum_{i=2}^{n} \frac{1}{i} \sum_{j=1}^{i} E[X_{i-j}]
     = \sum_{i=2}^{n} \frac{1}{i} \sum_{t=0}^{i-1} E[X_t],
\end{equation}
where $X_t := l.t.r.m.(\sigma)$ for a random $\sigma\in S_t$.

Recall from~\cite[Corollary 1.3.8]{ECI} that
$$%\begin{equation}\label{e.1.2}
\sum_{\sigma\in S_t} q^{l.t.r.m.(\sigma)}=\prod_{k=1}^{t} (q+k-1).
$$%\end{equation}
It follows that, for $t\ge 1$,
\begin{eqnarray*}%\label{e.3}
E[X_t] 
&=& \frac{1}{t!} \sum_{\sigma\in S_t} l.t.r.m.(\sigma)
 =  \frac{1}{t!} \left.\left(\frac{d}{dq}
    \sum_{\sigma\in S_t} q^{l.t.r.m.(\sigma)}\right)\right|_{q=1}\cr
&=& \frac{1}{t!} \left.\left(\frac{d}{dq}
    \prod_{k=1}^{t} (q+k-1)\right)\right|_{q=1}
 =  \frac{1}{t!} \sum_{r=1}^{t} \prod_{{1\le k\le t}\atop k\ne r} k 
 =  \sum\limits_{r=1}^{t} \frac{1}{r}.
\end{eqnarray*}
Of course, $E[X_0] = 0$.
Substituting these values into (\ref{e.2}) gives
$$
E[X] = \sum_{i=2}^{n} \sum_{t=1}^{i-1} \sum_{r=1}^{t} \frac{1}{i \cdot r}
$$
and this is equivalent (with $j=t+1$ and $k=r+1$) to the first formula
in the statement of the theorem. 

The second formula may be obtained through the following manipulations:
\begin{eqnarray*}
E[X] 
&=& \sum_{i=2}^{n} \sum_{j=2}^{i} \sum_{k=2}^{j} \frac{1}{i \cdot (k-1)}
 =  \sum_{2\le k\le j\le i\le n} \frac{1}{i \cdot (k-1)}\cr
&=& \sum_{2\le k\le i\le n} \frac{i-k+1}{i \cdot (k-1)}
 =  \sum_{2\le k\le i\le n} \left(\frac{1}{k-1} - \frac{1}{i}\right)\cr
&=& \sum_{2\le k\le n} \frac{n-k+1}{k-1} - \sum_{2\le i\le n} \frac{i-1}{i}\cr
&=& n\sum_{k=2}^{n} \frac{1}{k-1} - (n-1) - (n-1) + \sum_{i=2}^{n} \frac{1}{i}\cr
&=& n\sum_{i=1}^{n} \frac{1}{i} - 2n + \sum_{i=1}^{n} \frac{1}{i}.
\end{eqnarray*}

\qed

\noindent{\bf Proof of Corollary~\ref{t.e1}.} Notice that
$$
\sum\limits_{i=1}^n \frac{1}{i} = \ln n + O(1).
$$
(The next term in the asymptotic expansion is Euler's constant.)
Substitute into Theorem~\ref{t.e11} to obtain the desired result.

\qed

\section{Generalized Down Degrees}%Strong Descent Numbers

%\subsection{Definitions}

\begin{df}
For $\pi\in S_n$ and $1\le r< n$ let
$$
D_{-}^{(r)}(\pi):=\{t_{a,b}|\ \ell(\pi)>\ell(t_{a,b}\pi)>\ell(\pi)-2r\}
$$
the {\em $r$-th strong descent set} of $\pi$.

Define the {\em $r$-th
%strong descent number
down degree} as
$$
d_{-}^{(r)}(\pi):=\# D_{-}^{(r)}(\pi).
$$
%the cardinality of the $r$-th strong descent set.
\end{df}

\begin{exa}
The first strong descent set and
%number
down degree are those studied in the previous section; namely,
$D_{-}^{(1)}(\pi)=D_{-}(\pi)$ and $d_{-}^{(1)}(\pi)=d_{-}(\pi)$.

The $(n-1)$-th strong descent set is the set of {\em inversions}:
$$
D_{-}^{(n-1)}(\pi)=\{t_{a,b}\,|\, a<b,\,\pi^{-1}(a)>\pi^{-1}(b)\}.
$$
Thus
$$
d_{-}^{(n-1)}(\pi)=\inv(\pi),
$$
the inversion number of $\pi$.
\end{exa}

\begin{obs}\label{t.gob1}
For every $\pi\in S_n$ and $1\le a<b\le n$, $t_{a,b}\in
D_{-}^{(r)}(\pi)$ if and only if $\pi=[\dots,b,\dots,a,\dots]$ and
there are less than $r$ letters between the positions of $b$ and
$a$ in $\pi$ whose value is between $a$ and $b$.
\end{obs}

\begin{exa}
Let  $\pi=[7,9,5,2,3,8,4,1,6]$. Then
$$
D_{-}^{(1)}(\pi)=\{ t_{6,7}, t_{6,8}, t_{1,4}, t_{1,3}, t_{1,2}, t_{4,8}, t_{4,5},
t_{8,9}, t_{3,5}, t_{2,5}, t_{5,9}, t_{5,7} \}
$$
and
$$
D_{-}^{(2)}(\pi)=D_{-}^{(1)}(\pi)\cup \{ t_{6,9}, t_{1,8}, t_{4,9}, t_{4,7}, t_{3,9},
t_{3,7}, t_{2,9}, t_{2,7} \}.
$$
\end{exa}

\begin{cor}
For every $\pi\in S_n$ and $1\le r<n$
$$
d_{-}^{(r)}(\pi)=d_{-}^{(r)}(\pi^{-1}).
$$
\end{cor}

\noindent{\bf Proof.} By Observation~\ref{t.gob1}, $t_{a,b}\in
D_{-}^{(r)}(\pi)$ if and only if $t_{\pi^{-1}(a),\pi^{-1}(b)}\in
D_{-}^{(r)}(\pi^{-1})$.

\qed

%\begin{pro}\label{t.gdetermine}
%For every $r$, the $r$-th strong descent set of a permutation in
%$S_n$ uniquely determines the permutation.
%\end{pro}
%
%The proof is similar to the proof of
%Proposition~\ref{t.determine} and is omitted.

%By induction on $n$. Proposition clearly holds for $n=1$.
%
%Let $\pi$ be a permutation in $S_n$ and $\bpi$ be the
%permutation obtained by omitting $n$. By the induction hypothesis
%$\bpi$ is determined by the subset
%$$
%D_{-}^{(r)}(\bpi)=D_{-}^{(r)}(\pi)\setminus \{t_{a,n}| 1\le a<n\}.
%$$
% Hence it suffices to determine the position of $n$ in
%$\pi$.
%
%Now, if $j$ is the position of $n$ in $\pi$ then $\pi(j+1)$ is a
%the rightmost letter in $\pi$ among the letters in the set
%$$
% \{a|\ t_{a,n}\in D_{-}(\pi)\}.
%$$
%This determines the position of $n$ in $\pi$, completing the
%proof.
%
%\qed

%\subsection{Maximum}

\begin{df}
The {\em $r$-th strong descent graph} of $\pi\in S_n$, denoted $\Gamma_{-}^{(r)}(\pi)$,
is the graph whose set of vertices is $\{1,\dots,n\}$ and whose
set of edges is
$$
\{\{a,b\}|\ t_{a,b}\in D_{-}^{(r)}(\pi)\}.
$$
\end{df}

The following lemma generalizes Lemma~\ref{t.triangle}.

\begin{lem}\label{t.gtriangle}
For every $\pi\in S_n$, the graph $\Gamma_{-}^{(r)}(\pi)$
contains no subgraph isomorphic to the complete graph $K_{r+2}$.
\end{lem}

\noindent{\bf Proof.} 
Assume that there is a subgraph in
$\Gamma_{-}^{(r)}(\pi)$ isomorphic to $K_{r+2}$.
Then there exists a decreasing subsequence
$$
n\ge a_1>a_2>\cdots>a_{r+2}\ge 1
$$
such that for all $1\le i<j\le r+2$, $t_{a_i,a_j}$ are $r$-th
strong descents of $\pi$. In particular, for every $1\le i<r+2$,
$t_{a_i,a_{i+1}}$ are $r$-th strong descents of $\pi$. This
implies that, for every $1\le i<r+2$, $a_{i+1}$ appears to the right of
$a_i$ in $\pi$. Then, by Observation~\ref{t.gob1}, $t_{a_1,a_{r+2}}$
is not an $r$-th strong descent. Contradiction.

\qed

\begin{cor}\label{t.gmax}
For every $1\le r< n$,
$$
max\{d_{-}^{(r)}(\pi) \,|\, \pi\in S_n\}\le t_{r+1}(n)\le 
{r+1\choose 2}\left(\frac{n}{r+1}\right)^2 .
$$
\end{cor}

\noindent{\bf Proof.} 
Combining Tur\'an's Theorem together with Lemma~\ref{t.gtriangle}.

\qed

\medskip

Note that for $r=1$ and $r=n-1$ equality holds in
Corollary~\ref{t.gmax}.

%\begin{prb}
%For which $r$-s equality holds in Corollary~\ref{t.gmax} ?
%\end{prb}

%\subsection{Expectation and Distribution}

%\section{Final Remarks}

\begin{rem} For every $\pi\in S_n$ let $\bpi$ be the permutation
obtained from $\pi$ by omitting the value $n$. 
If $j$ is the position of $n$ in $\pi$ then
$$
d_{-}^{(r)}(\pi)-d_{-}^{(r)}(\bpi)
$$
equals the number of $(r-1)$-th almost left-to-right minima in the $(j-1)$-th
suffix of $\bpi$, see e.g.~\cite{RR}. This observation may be applied to
calculate the expectation of $d_{-}^{(r)}(\pi)$.
\end{rem}

%%{\bf APPLY THIS OBSERVATION TOGETHER WITH RESULTS ON ALMOST
%%L.T.R.MAX TO CALCULATE EXPECTATION OF $d_{-}^{(r)}(\pi)$ ETC}
%
%%\section{Other Types}
%
%\begin{rem}
%By results of Incitti~\cite{Incitti}, Propositions~\ref{t.determine}
%and~\ref{t.gdetermine} hold for type $B$.
%\end{rem}
%%ARE THEY CORRECT FOR GENERAL TYPE ??}

\bigskip

\noindent{\bf Acknowledgements.} 
The concept of strong descent graph came up during conversations with
Francesco Brenti. Its name and certain other improvements were suggested
by Christian Krattenthaler. Thanks also to Nathan Reading, Amitai Regev,
Alexander Yong, and the anonymous referees.

\end{document}